\documentclass[11pt]{article}

\usepackage{graphics}
\usepackage{amsfonts}
\usepackage{amssymb}
\usepackage{amsmath}
\usepackage{amsthm}

\usepackage{url}
\usepackage{color}
\usepackage{multirow}
\usepackage{enumerate}
\usepackage{hyperref}

\newtheorem{thm}{Theorem}
\newtheorem{lem}{Lemma}
\newtheorem{cor}{Corollary}

\newtheorem{deff}{Definition}

\newcommand{\R}{\mathbb{R}}

\newcommand{\cA}{\mathcal{A}}

\DeclareMathOperator{\tr}{tr}
\DeclareMathOperator{\diag}{diag}
\DeclareMathOperator{\Diag}{Diag}
\DeclareMathOperator{\aut}{aut}
\DeclareMathOperator{\spen}{span}

\renewcommand{\binom}[2]{\genfrac(){0pt}{}{#1}{#2}}

\title{Semidefinite programming and eigenvalue bounds for the graph partition problem\footnote{This version is published in
Mathematical Programming, DOI:\href{http://dx.doi.org/10.1007/s10107-014-0817-6}{10.1007/s10107-014-0817-6}}}
\author{Edwin R. van Dam \thanks{Department of Econometrics and OR, Tilburg
University, The Netherlands. {\tt edwin.vandam@uvt.nl} }
 \and {Renata Sotirov}\thanks{Department of Econometrics and OR, Tilburg
University, The Netherlands. {\tt r.sotirov@uvt.nl} }}

\date{}
\begin{document}
\maketitle

\begin{abstract}
The graph partition problem is the problem of partitioning the vertex set of a graph into a fixed number of
sets of given sizes such that the sum of weights of edges joining different sets is optimized.
In this paper we
simplify a known matrix-lifting semidefinite programming relaxation of the graph partition problem for
several classes of graphs and also show how to aggregate additional triangle and independent set constraints
for graphs with symmetry. We present  an eigenvalue bound for the graph partition problem of a strongly
regular graph, extending a similar result for the equipartition problem.
We also derive a linear programming bound of the graph partition problem
for certain Johnson and Kneser graphs. Using what we call the Laplacian algebra of a graph, we derive an
eigenvalue bound for the graph partition problem that is the first known closed form bound that is applicable
to any graph, thereby extending a well-known result in spectral graph theory. Finally, we strengthen a known
semidefinite programming relaxation of a specific quadratic assignment problem and the
above-mentioned matrix-lifting semidefinite programming relaxation by adding two constraints that correspond
to assigning two vertices of the graph to different parts of the partition. This strengthening performs well
on highly symmetric graphs when other relaxations provide weak or trivial bounds.
\end{abstract}


\section{Introduction}\label{sec:intro}

The graph partition  problem (GPP) is the problem of partitioning the vertex set of a graph into a fixed
number, say $k$, of sets of given sizes such that the sum of weights of edges joining different sets is
optimized. Here we also refer to the described GPP problem as the $k$-partition problem. The GPP is a NP-hard
combinatorial optimization problem, see \cite{GarJoSt:76}. It has many applications such as VLSI design
\cite{Len:90}, parallel computing \cite{BisHenKa:00,HeKo:00,Simon91}, network partitioning
\cite{FidMatt:82,Sanch89}, and floor planing \cite{DaKuh87}. For recent advances in graph
partitioning, we refer to \cite{BMSSS:13}.

There are several approaches for deriving bounds for the GPP. Here we are interested in eigenvalue and
semidefinite programming (SDP) bounds. Donath and Hoffman  \cite{DoHo:73} derived an eigenvalue-based
bound for the GPP that was further improved by Rendl and Wolkowicz \cite{ReWo:95}. Falkner, Rendl, and
Wolkowicz \cite{FeReWo:92} derived a closed from bound for the minimum $k$-partition problem when
$k=2,3$ by using the bound from \cite{ReWo:95}. Their bound for $k=2$ coincides with a well-established
result in spectral graph theory; see e.g., Juvan and Mohar \cite{JuMo:92}. Alizadeh \cite{Aliz95} and Karish
and Rendl \cite{KarRend:98} showed that the Donath-Hoffman bound from \cite{DoHo:73} and the Rendl-Wolkowicz
bound from \cite{ReWo:95}, respectively, can be reformulated as semidefinite programs. Other SDP relaxations
of the GPP were derived in \cite{KarRend:98,WoZh:99,KaReCl:00,dKPaDoSo:10,Sot11}. For a comparison of all
these relaxations, see \cite{Sot10,Sot11}.

Armbruster, Helmberg, F\"ugenschuh, and Martin \cite{ArmHeFueMa:11}  evaluated the strength of a
branch-and-cut framework for linear and semidefinite relaxations of the minimum graph bisection problem on
large and sparse instances. Their results show that in the majority of the cases the semidefinite approach
outperforms the linear one. This is very encouraging since SDP relaxations are widely
believed to be of use only for instances that are small and dense.
The aim of this paper is to further investigate eigenvalue and SDP bounds for the GPP. \\

\noindent {\bf Main results and outline}

\vspace{0.2cm} \noindent Symmetry in graphs is typically ascribed to symmetry coming from the automorphism
group of the graph. However, symmetry may also be interpreted (and exploited) more broadly as what we call
combinatorial symmetry. In Section \ref{sec:SyminGr} we explain both types of symmetry. In particular, in
Section \ref{sec:perm} we explain symmetry coming from groups, while in Section \ref{sec:comb} we describe
combinatorial symmetry and related coherent configurations. The associated coherent algebra is an algebra
that can be exploited in many combinatorial optimization problems, in particular the GPP. In the case that a
graph has no (or little) symmetry, one can still exploit the algebraic properties of its Laplacian matrix.
For this purpose, we introduce the Laplacian algebra and list its properties in Section \ref{Ex:LaplAlg}.

In Section \ref{sect:EffRelax} we simplify the matrix-lifting SDP relaxation from \cite{Sot11} for different
classes of graphs and also show how to aggregate triangle and independent set constraints when possible,
see Section \ref{sec:aggregateConstr}. This approach enables us, for example, to solve the SDP relaxation
from \cite{Sot11} with an additional $3\binom{n}{3}$ triangle constraints in less than a second (!) for
highly symmetric graphs with $n=100$ vertices. In Section \ref{sec:SRG}  we present  an eigenvalue bound
 for the GPP of a strongly regular graph (SRG). This result is an extension of the result by De
Klerk et al.~\cite{{dKPaDoSo:10},deKlSotNaTr:10} where an eigenvalue bound is derived for the graph
equipartition problem for a SRG. In Section \ref{sec:SRG} we also show that for all
SRGs except for the pentagon, the bound from \cite{Sot11} does not improve by adding triangle inequalities.
In Section \ref{sec:JohnKne} we derive a linear program  that is equivalent to the SDP relaxation
from \cite{Sot11} when the graph under consideration is a Johnson or Kneser graph on triples.

In Section \ref{sec:AnyGaph} we derive an eigenvalue bound for the GPP for any, not
necessarily highly symmetric, graph. This is the first known closed form bound for the minimum $k$-partition
when $k>3$ and for the maximum  $k$-partition when  $k> 2$ that is applicable to any graph. Our result is a
generalization of a well-known result in spectral graph theory for the $2$-partition problem to any $k$-partition problem.

In Section \ref{sec:impVec},  we derive a new  SDP relaxation  for the GPP that is suitable for graphs with
symmetry. The new relaxation is a strengthened SDP relaxation of a specific quadratic assignment
problem (QAP) by Zhao, Karisch, Rendl, and Wolkowicz \cite{ZhKaReWo:98} by adding two constraints that
correspond to assigning two vertices of the graph to different parts of the partition. The new bound performs
well on highly symmetric graphs when other SDP relaxations provide weak or trivial bounds. This is probably
due to the fact that fixing breaks (some) symmetry in the graph under consideration. Finally, in Section
\ref{sec:impMat} we show how to strengthen the matrix-lifting relaxation from \cite{Sot11}  by adding a
constraint that corresponds to assigning two vertices of the graph to different parts of the partition. The
new matrix-lifting SDP relaxation is not dominated by the relaxation from \cite{ZhKaReWo:98}, or vice versa.

The numerical results in Section \ref{sec:NumRez} present the high potential of the new
bounds. None of the presented SDP bounds strictly dominates any of the other bounds for all
tested instances. The results indicate that breaking symmetry strengthens the bounds from
\cite{Sot11,ZhKaReWo:98} when the triangle and/or independent set constraints do not (or only slightly)
improve the bound from \cite{Sot11}. For the cases that the triangle and/or independent set constraints
significantly improve the bound from \cite{Sot11}, the fixing approach does not seem to be very effective.

\section{Symmetry and matrix algebras} \label{sec:SyminGr}

A matrix $*$-algebra is a set of matrices that is closed under addition, scalar multiplication, matrix
multiplication, and taking conjugate transposes. In \cite{GatPa:04,EdeKl:10,deKlDoPa:11} and others, it was
proven that  one can restrict optimization of an SDP problem  to feasible points in a matrix $*$-algebra that
contains the data matrices of that problem. In particular, the following theorem is proven.
\begin{thm} {\em \cite{deKlDoPa:11}} \label{thmAlg}
Let $\mathcal A$ denote a matrix $*$-algebra that contains the data matrices of an SDP problem as well as the
identity matrix. If the SDP problem has an optimal solution, then it has an optimal solution in $\mathcal A$.
\end{thm}
When the matrix $*$-algebra has small dimension, then one can exploit a basis of the algebra to reduce the
size of the SDP considerably, see e.g, \cite{EdeKl:10,dKPaDoSo:10,KOP}. In the recent papers
\cite{EdeKl:10,dKPaDoSo:10,deKlSot:10,deKlSotNaTr:10}, the authors considered matrix $*$-algebras that
consist of matrices that commute with a given set of permutation matrices that correspond to  automorphisms.
Those $*$-algebras have a basis of
$0$-$1$ matrices that can be efficiently computed. However, there exist also such $*$-algebras that are not
coming from permutation groups, but from the `combinatorial symmetry', as we shall see below. We also
introduce the Laplacian algebra in order to obtain an eigenvalue bound that is suitable for any graph.

Every matrix $*$-algebra $\cA$ has a canonical block-diagonal structure. This is a consequence of the theorem by Wedderburn
\cite{Wedderburn} that states that there is a $*$-isomorphism
\[
\varphi:\cA \longrightarrow \oplus_{i=1}^p  \mathbb{C}^{n_i\times n_i},
\]
i.e., a bijective linear map that preserves multiplication and conjugate transposition.
One can  exploit such a $*$-isomorphism in order to further reduce the size of an SDP.

\subsection{Symmetry from automorphisms } \label{sec:perm}
An {\em automorphism} of a graph $G=(V,E)$ is a bijection $\pi:V \rightarrow V$ that preserves edges, that
is, such that $\{\pi(x),\pi(y)\} \in E$ if and only if $\{x,y\}\in E$. The set of all automorphisms of $G$
forms a group under composition; this is called the {\em automorphism group} of $G$. The {\em orbits} of the
action of the automorphism group acting on $V$ partition the vertex set $V$; two vertices are in the same
orbit if and only if there is an automorphism mapping one to the other. The graph $G$ is {\em
vertex-transitive} if its automorphism group acts transitively on vertices, that is, if for every two
vertices, there is an automorphism that maps one to the other (and so there is just one orbit of vertices).
Similarly, $G$ is {\em edge-transitive} if its automorphism group acts transitively on edges. Here, we
identify the automorphism group of the graph with the automorphism group of its adjacency matrix. Therefore,
if $G$ has adjacency matrix $A$ we will also refer to the automorphism group of the graph as $\aut(A) :=\{P
\in \Pi_n: AP=PA\}$, where $\Pi_n$ is the set of permutation matrices of size $n$.

Assume that $\mathcal G$ is a subgroup of the automorphism group of $A$. Then the centralizer ring (or
commutant) of ${\mathcal G}$, i.e.,
${\mathcal A}_{\mathcal G}=\{ X\in \R^{n\times n}: XP=PX, ~\forall P \in {\mathcal G} \}$
is a matrix $*$-algebra that contains $A$. One may obtain a basis for ${\mathcal A}_{\mathcal G}$
from the orbitals (i.e., the orbits of the action of $\mathcal G$ on ordered pairs of vertices) of the group $\mathcal G$.
This basis, say $\{A_1$,\ldots, $A_r\}$ forms a so-called coherent configuration.
\begin{deff}[Coherent configuration]
\label{def:coherent config} A set of zero-one $n\times n$ matrices $ \{A_1,\ldots, A_r\}$ is called a
\emph{coherent configuration} of rank $r$ if it satisfies the following properties:
\begin{enumerate}[(i)]
\item $\sum_{i \in \mathcal{I}} A_i = I$ for some index set $\mathcal{I} \subset \{1,\ldots,r\}$ and
    $\sum_{i=1}^r A_i = J$,
\item $A_i^\mathrm{T} \in \{A_1,\ldots, A_r\}$ for $i=1,\ldots, r$,
\item There exist $p^h_{ij}$, such that   $A_iA_j =\sum_{h=1}^r     p^h_{ij}A_h$ for $i,j\in\{1,\ldots,r\}$.
\end{enumerate}
\end{deff}
As usual, the matrices $I$ and $J$ here denote the identity matrix and all-ones matrix, respectively.
 We call $\mathcal{A}:=\spen\{A_1,\dots,A_r\}$ the associated {\em coherent algebra}, and this is clearly a
matrix $*$-algebra. Note that in the case that  $A_1,\dots,A_r$ are derived as orbitals of the group
$\mathcal G$, it follows indeed that $\mathcal{A}={\mathcal A}_{\mathcal G}$. If the coherent configuration
is commutative, that is, $A_iA_j=A_jA_i$ for all $i,j=1,\dots,r$, then we call it a (commutative) {\it
association scheme}. In this case, $\mathcal{I}$ contains only one index, and it is common to call this index
$0$ (so $A_0=I$), and $d:=r-1$ the number of classes of the association scheme.

In the case of an association scheme, all matrices can be diagonalized simultaneously, and the corresponding
$*$-algebra has a canonical diagonal structure $\oplus_{j=0}^d \mathbb{C}$. The $*$-isomorphism $\varphi$ is then given
by $\varphi(A_i)=\oplus_{j=0}^d P_{ji}$, where $P_{ji}$ is the eigenvalue of $A_i$ on the $j$-th eigenspace. The matrix
$P=(P_{ji})$ of eigenvalues is called the {\em eigenmatrix} or {\em character table} of the association scheme.

Centralizer rings are typical examples of coherent algebras, but not the only ones. In
general, the centralizer ring of the automorphism group of $A$ is not the smallest coherent algebra
containing $A$, even though this is the case for well-known graphs such as the  Johnson and Kneser graphs
that we will encounter later in this paper. We could say that, in general, the smallest coherent
configuration captures more symmetry than that coming from automorphisms of the graph. In this case, we say
that there is more combinatorial symmetry.

\subsection{Combinatorial symmetry} \label{sec:comb}

Let us look at coherent configurations and the combinatorial symmetry that they capture in more detail. One
should think of the (non-diagonal) matrices $A_i$ of a coherent configuration as the adjacency matrices of
(possibly directed) graphs on $n$ vertices. The diagonal matrices represent the different `kinds' of vertices
(so there are $|\mathcal{I}|$ kinds of vertices; these generalize the orbits of vertices under the action of
the automorphism group). The non-diagonal matrices $A_i$ represent the different `kinds' of edges and
non-edges.

In order to identify the `combinatorial symmetry' in a graph, one has to find a coherent configuration
(preferably of smallest rank) such that the adjacency matrix of the graph is in the corresponding coherent
algebra $\mathcal{A}$. As mentioned before, not every coherent configuration comes from the orbitals of a
permutation group. Most strongly regular graphs --- a small example being the  Shrikhande graph
--- indeed give rise to such examples. A (simple, undirected, and loopless) $\kappa$-regular graph $G=(V,E)$
on $n$ vertices is called {\em strongly regular} with parameters $(n, \kappa, \lambda, \mu)$ whenever it is
not complete or edgeless and every two distinct vertices have $\lambda$ or $\mu$ common neighbors, depending
on whether the two vertices are adjacent or not, respectively. If $A$ is the adjacency matrix of $G$, then
this definition implies that $A^2=\kappa I+\lambda A+\mu(J-I-A)$, which implies furthermore that
$\{I,A,J-I-A\}$ is an association scheme. The combinatorial symmetry thus tells us that there is one kind of
vertex, one kind of edge, and one kind of non-edge. For the Shrikhande graph, a strongly regular graph with
parameters $(16,6,2,2)$ (defined by $V=\mathbb{Z}_4^2$, where two vertices are adjacent if their difference
is $\pm(1,0),\pm(0,1),$ or $\pm(1,1)$) however, the automorphism group indicates that there are two kinds of
non-edges (depending on whether the two common neighbors of a non-edge are adjacent or not), and in total
there are four (not three) orbitals. Doob graphs are direct products of $K_4$s and Shrikhande graphs, thus
generalizing the Shrikhande graph to association schemes with more classes. In many optimization problems the
combinatorial symmetry, captured by the concept of a coherent configuration or association scheme, can be
exploited, see Section \ref{sec:aggregateConstr} and e.g.,~\cite{KOP,{GoeRend:99}}.
In Section \ref{ces:CombSymNum}, we will mention some numerical results for graphs that have
more combinatorial symmetry than symmetry coming from automorphisms.

\subsection{The Laplacian algebra} \label{Ex:LaplAlg}

Let $A$ be an adjacency matrix of a connected graph $G$ and $L:=\Diag(Au_n)-A$  the Laplacian matrix of the graph.
We introduce the matrix $*$-algebra consisting of all polynomials in $L$, and call this algebra the Laplacian algebra
$\mathcal L$. This algebra has a convenient basis of idempotent matrices that are formed from an orthonormal basis of
eigenvectors corresponding to the eigenvalues of $L$. In particular, if the distinct eigenvalues of $L$ are denoted by
$0=\lambda_0 < \lambda_1 < \cdots < \lambda_d$, then we let $F_i=U_iU_i^{\mathrm{T}}$, where $U_i$ is a matrix
having as columns an orthonormal basis of the eigenspace of $\lambda_i$, for $i=0, \ldots,d$.
Then $\{F_0, \dots, F_d\}$ is a basis of $\mathcal L$ that satisfies the following properties:
\begin{itemize}
\item $F_0= \frac{1}{n} J$,  ~$\sum\limits_{i=0}^dF_i=I$, $\sum\limits_{i=0}^d \lambda_iF_i~=L$
\item $F_iF_j = \delta_{ij} F_i$, ~$\forall i,j$
\item $F_i=F_i^*$, ~$\forall i$.
\end{itemize}
Note that  $\tr F_i=f_i$, the multiplicity of eigenvalue $\lambda_i$ of $L$, for all $i$. Clearly,
the operator $P$, where
$$P(Y)=\sum_{i=0}^d \frac{\tr YF_i}{f_i} F_i$$
is the orthogonal projection onto $\mathcal L$.

We note that the Laplacian algebra of a strongly regular graph is the same as the corresponding coherent algebra $\spen
\{I,A,J-I-A\}$  (and a similar identity holds for graphs in association schemes).

\section{The graph partition problem}\label{sec:gpp}

The  minimum  (resp.~maximum) graph partition problem  may be formulated as follows.
Let $G=(V,E)$ be an undirected graph  with vertex set $V$, where $|V|=n$ and edge set $E$, and $k\geq 2$ be a given integer.
The goal is to find a partition  of the vertex set into $k$ (disjoint) subsets $S_1,\ldots, S_k$ of
specified sizes $m_1\geq \ldots \geq m_k $, where $\sum_{j=1}^k m_j =n$, such that the sum of weights of edges joining different sets
$S_j$ is minimized (resp.~maximized).
The case when $k=2$ is known  as the {\em graph bisection problem} (GBP).
If all $m_j$ ($j=1,\ldots,k$) are equal, then we refer to the associated problem  as the {\em graph equipartition problem} (GEP).

We denote by $A$ the adjacency matrix of $G$. For a given partition of the graph into $k$ subsets,   let   $X=(x_{ij})$ be
the $n\times k$ matrix defined by
\[
x_{ij} = \left \{
\begin{array}{ll}
1 &  \mbox{if $i\in S_j$ } \\
0 & \mbox{otherwise}.
\end{array} \right .
\]
Note that the $j$th column of $X$ is the characteristic vector of $S_j$. The sum of weights of edges joining
different sets, i.e., the cut of the partition, is equal to $\frac{1}{2} \tr A(J_n -XX^{\mathrm{T}})$. Thus,
the minimum  GPP problem can be formulated as
\begin{equation} \label{GPP}
\begin{array}{ll}
  \min  &  \frac{1}{2} \tr A(J_n -XX^{\mathrm{T}}) \\[1ex]
 {\rm s.t.} & Xu_k=u_n\\[1ex]
 & X^{\mathrm{T}}u_n=m\\[1ex]
& x_{ij}\in \{0,1\},~~\forall i,j,
\end{array}
\end{equation}
where $m=(m_1,\ldots, m_k)^\mathrm{T}$, and $u_k$ and $u_n$ denote all-ones vectors of sizes $k$ and $n$, respectively.
It is easy to show that if $X$ is feasible for (\ref{GPP}), then
\begin{equation} \label{obj1}
\frac{1}{2}\tr A(J_n -XX^{\mathrm{T}})=\frac{1}{2} \tr LXX^{\mathrm{T}},
\end{equation}
where  $L$ is the Laplacian matrix of the graph. We will use this alternative expression for the objective in Section \ref{sec:AnyGaph}.

\section{A simplified and improved SDP relaxation for the GPP} \label{sect:EffRelax}

In \cite{Sot11}, the second author derived a matrix lifting SDP relaxation for the GPP.
Extensive numerical results in \cite{Sot11} show that the matrix lifting SDP relaxation for the GPP provides competitive bounds and
is solved significantly faster than any other known SDP bound for the  GPP.
The goal of this section is to further simplify the mentioned relaxation for highly symmetric graphs.
Further, we show here how to aggregate, when possible, certain types of (additional) inequalities to obtain stronger bounds.

The matrix lifting relaxation in \cite{Sot11} is obtained after linearizing the objective function
$\tr A(J_n -XX^{\mathrm{T}})$ by replacing $XX^{\mathrm{T}}$ with a new variable $Y$,  and approximating the set
\[
{\rm conv} \left \{ XX^{\mathrm{T}}: X\in \R^{n\times k}, ~Xu_k=u_n, ~X^{\mathrm{T}}u_n=m, ~x_{ij}\in \{0,1\} \right \}.
\]
The following SDP relaxation for the GPP is thus obtained.
\[({\rm GPP}_{\rm m})\quad
\begin{array}{rl}
\min & \frac{1}{2} \tr A(J_n -Y) \\[1ex]
{\rm s.t.} & \diag(Y)=u_n \\
& \tr JY =\sum\limits_{i=1}^k m_i^2  \\
& kY - J_n\succeq 0, ~~Y\geq 0.
\end{array}
\]
We observe the following simplification of the relaxation ${\rm GPP}_{\rm m}$
for the bisection problem.
\begin{lem}
For the case of the  bisection problem the nonnegativity constraint on the matrix variable in ${\rm GPP}_{\rm m}$ is redundant.
\end{lem}

\noindent
{\em Proof}.  Let $Y$ be feasible for ${\rm GPP}_{\rm m}$ with $k=2$.
We define $Z:= 2Y - J_n$. Now from $\diag(Y)=u_n$ it follows that $\diag(Z)=u_n$.
Because $Z \succeq 0$, it follows that $-1\leq z_{ij} \leq 1$, which implies indeed that $y_{ij} \geq 0$. \qed\\

In order to strengthen ${\rm GPP}_{\rm m}$, one can add the triangle constraints
\begin{equation} \label{triangle}
y_{ab} + y_{ac}\leq 1 + y_{bc}, \quad \forall (a,b,c).
\end{equation}
For a given triple $(a,b,c)$ of (distinct) vertices, the constraint \eqref{triangle} ensures that if $a$ and $b$
belong to the same set of the partition and so do $a$ and $c$, then also $b$ and $c$ do so. There are  $3\binom{n}{3}$ inequalities of type \eqref{triangle}.
For future reference, we refer to ${\rm GPP}_{\rm m\triangle}$ as the SDP relaxation that is obtained
from ${\rm GPP}_{\rm m}$ by adding the triangle constraints.

One can also add  to ${\rm GPP}_{\rm m}$  and/or ${\rm GPP}_{\rm m\triangle}$ the independent set constraints
\begin{equation} \label{indepSet}
\sum\limits_{a<b, ~a,b\in W} y_{ab}\geq 1, ~\mbox{for all}~ W ~\mbox{with}~ |W|=k+1.
\end{equation}
These constraints ensure that the graph with adjacency matrix $Y$ has no independent set ($W$) of size $k+1$.
There are $\binom{n}{k+1}$ inequalities of type \eqref{indepSet}. For future reference, we refer to ${\rm GPP}_{\rm m-ind}$ as the SDP relaxation that is obtained from ${\rm GPP}_{\rm m}$ by adding  the independent set constraints.

Constraints \eqref{triangle} and  \eqref{indepSet} are also  used by Karish and Rendl \cite{KarRend:98} to strengthen
the SDP relaxation for the graph equipartition problem, and by the second author
 \cite{Sot11} to strengthen the SDP relaxation for the (general) graph partition problem.
By adding constraints \eqref{triangle} and/or \eqref{indepSet} to ${\rm GPP}_{\rm m}$, one obtains --- in general ---
stronger relaxations that are computationally more demanding than ${\rm GPP}_{\rm m}$.
In the following sections we will show how to efficiently compute, for graphs with symmetry, all above derived relaxations.

\subsection{Symmetry and aggregating triangle and independent set constraints} \label{sec:aggregateConstr}

It is well known how to exploit the symmetry in problems such as ${\rm GPP}_{\rm m}$ by using coherent configurations (or association schemes). Aggregating triangle inequalities was suggested by Goemans and  Rendl \cite{GoeRend:99} in the context of the maximum cut problem for graphs in association schemes.  Surprisingly, the suggestion by Goemans and Rendl was not followed so far in the literature, as far as we know. Here we will extend the approach successfully to coherent configurations. Moreover, we will aggregate the independent set inequalities \eqref{indepSet} for the case $k=2$.

Let us now consider graphs with symmetry, and assume that the data matrices of ${\rm GPP}_{\rm m}$ belong to the coherent algebra of a coherent configuration $\{ A_1,\ldots, A_r\}$. We will first show how this allows us to efficiently solve ${\rm GPP}_{\rm m}$, and subsequently how to aggregate additional triangle and/or
independent set constraints.

Because of our assumption, we may consider $Y=\sum_{j=1}^r y_jA_j$ (see  Theorem \ref{thmAlg}) and the SDP relaxation ${\rm GPP}_{\rm m}$  reduces to
\begin{equation} \label{RSred}
\begin{array}{rl}
\min & \frac{1}{2} \tr AJ_n - \frac{1}{2}\sum\limits_{j=1}^r y_j \tr AA_j \\[1.5ex]
{\rm s.t.} & \sum\limits_{j\in {\mathcal I}} y_j \diag(A_j) =u_n \\[1ex]
& \sum\limits_{j=1}^r y_j \tr JA_j = \sum\limits_{i=1}^k m_i^2  \\[2.5ex]
& k\sum\limits_{j=1}^r y_jA_j - J_n\succeq 0, ~~y_j\geq 0,~~ j=1,\ldots,r,
\end{array}
\end{equation}
where $\mathcal I$ is the subset of $\{1,\ldots, r\}$ that contains elements of the coherent configuration with nonzero diagonal
(as in Definition \ref{def:coherent config}).
Note that \eqref{RSred}  solves significantly faster than ${\rm GPP}_{\rm m}$  when $r \ll n^2/2$.
Also, the linear matrix inequality in  \eqref{RSred} can be block-diagonalized.
In the following sections we will show that the SDP relaxation \eqref{RSred} can be further simplified for some special types of graphs.

Next, we will reduce ${\rm GPP}_{\rm m\triangle}$ by adding aggregated triangle inequalities  to \eqref{RSred}. Because we cannot express a single triangle inequality in terms of the new variables in \eqref{RSred}, we consider  all inequalities, at once, of the same `type', as follows.
For a  given triple of distinct vertices $(a,b,c)$ consider the triangle inequality
$y_{ab} + y_{ac}\leq 1 + y_{bc}$.
If $(A_i)_{ab}=1$, $(A_h)_{ac}=1$, and $(A_j)_{bc}=1$, then we say that this triangle inequality is of type $(i,j,h)$
($i,j,h \in \{1,\ldots, r\} \setminus {\mathcal I}$; note that the indices $i,j,h$ are not necessary distinct).
From  Definition \ref{def:coherent config} (iii) it follows that if $(A_i)_{ab}=1$,
then the number of (directed) triangles containing the (directed) edge $(a,b)$ and for which  $(A_h)_{ac}=1$ and $(A_{j'})_{cb}$ is equal to $p^i_{hj'}$ (here $j'$ is the index for which  $A_{j'}=A_j^{\mathrm T}$).
Now, if $Y$ is feasible for \eqref{RSred}, then by summing all triangle inequalities of a given type $(i,j,h)$, the aggregated triangle inequality becomes
\begin{equation}\label{CohTriangAgg}
p^i_{hj'} \tr A_iY + p^h_{ij} \tr A_hY  \leq p^i_{hj'} \tr A_iJ + p^j_{i'h} \tr A_{j}Y.
\end{equation}
After exploiting the fact that $Y=\sum_{j=1}^r y_jA_j$,  the aggregated inequality \eqref{CohTriangAgg} reduces to a linear inequality that can be added to the relaxation
\eqref{RSred}. The number of aggregated triangle inequalities is bounded by $r^3$ which may be significantly smaller than $3\binom{n}{3}$.
So the SDP relaxation ${\rm GPP}_{\rm m\triangle}$ can be efficiently computed for small $r$.

For the bisection problem (i.e., $k=2$), the independent set constraints \eqref{indepSet} can be aggregated in a similar way as the triangle inequalities, and we obtain that
\begin{equation}\label{CohIndependAgg}
p^i_{hj'} \tr A_iY + p^h_{ij} \tr A_hY  + p^j_{i'h} \tr A_{j}Y \geq p^i_{hj'} \tr A_iJ.
\end{equation}
It is not clear how to aggregate the independent set constraints for $k\geq 3$. Note that in the case that the considered coherent configuration is an association scheme, all matrices are symmetric which simplifies the above aggregation processes.

\subsection{Strongly regular graphs} \label{sec:SRG}

In this section we derive a closed form expression for the  optimal objective value of the SDP relaxation ${\rm GPP}_{\rm m}$, for a strongly regular graph.
A similar approach was used in \cite{dKPaDoSo:10,{deKlSotNaTr:10}} to derive an eigenvalue bound for the equipartition problem from the SDP
relaxation presented by Karish and Rendl \cite{KarRend:98}.
Furthermore, we show that the triangle inequalities are redundant in ${\rm GPP}_{\rm m}$
 for connected SRGs, except for the pentagon.

Let $A$ be the adjacency matrix of a strongly regular graph $G$ with parameters $(n, \kappa, \lambda, \mu)$, see Section \ref{sec:comb}.
Using the matrix equation $A^2=\kappa I+\lambda A+\mu(J-I-A)$, we can determine the eigenvalues of the matrix $A$ from the
parameters of $G$, see e.g., \cite{BrHa}. Since $G$ is regular with valency $\kappa$, it follows that $\kappa$ is an
eigenvalue of $A$ with eigenvector $u_n$. The matrix $A$ has exactly two distinct eigenvalues associated with
eigenvectors orthogonal to $u_n$. These two eigenvalues are known as {\em restricted eigenvalues} and are usually
denoted by $r \geq 0$ and $s<0$. The character table of the corresponding association scheme is
\begin{equation} \label{charTabSRG}
P= \left (
\begin{array}{ccc}
1 & \kappa & n-1-\kappa \\[1ex]
1 & r & -1-r\\[1ex]
1 & s & -1-s
\end{array}
\right ).
\end{equation}
From Theorem \ref{thmAlg} it follows that there exists an optimal solution $Y$ to ${\rm GPP}_{\rm m}$ in the
coherent algebra spanned by $\{I, A, J-A-I \}$. Because of the constraints $\diag(Y)=u_n$ and $Y \geq 0$, there exist   $y_1, y_2 \geq 0$
such that
\begin{equation}\label{strY}
Y= I + y_1A + y_2(J-A-I),
\end{equation}
which we shall use to get an even simpler form than \eqref{RSred}.
The constraint $\tr JY = \sum_{i=1}^k m_i^2$ reduces to
\begin{equation} \label{eqtr}
n + n\kappa y_1 + (n^2-n\kappa -n)y_2 = \sum\limits_{i=1}^k m_i^2.
\end{equation}
Since the matrices $\{I, A, J-A-I \}$  may be simultaneously diagonalized,
the constraint $kY - J_n\succeq 0$ becomes a system of linear inequalities in the variables $y_1$ and $y_2$.
In particular, after exploiting \eqref{charTabSRG}, the constraint  $kY - J_n\succeq 0$ reduces to the three constraints
\begin{align}
k+k \kappa y_1+ k(n-\kappa-1)y_2 -n &\geq 0,  \label{lineq1} \\[1ex]
1+r y_1 - (r+1)y_2 &\geq 0,  \label{lineq2} \\[1ex]
1+s y_1-(s+1)y_2 &\geq 0. \label{lineq3}
\end{align}
Because $\sum\limits_{i=1}^k m_i^2 \geq n^2/k$ by Cauchy's inequality, \eqref{lineq1} is actually implied  by \eqref{eqtr}, so we may remove this first constraint.
It remains only to rewrite the objective function, i.e.,
\[
\frac{1}{2} \tr{A(J_n -Y)}= \frac{\kappa n(1-y_1)}{2}.
\]
To summarize, the SDP   bound ${\rm GPP}_{\rm m}$ can be obtained by solving the following linear
programming (LP) problem
\begin{equation} \label{LPsrg}
\begin{array}{ll}
\min &   \frac{1}{2}\kappa n(1-y_1)  \\[1ex]
{\rm s.t.} &  \kappa  y_1 + (n-\kappa -1)y_2 = \frac{1}{n}\sum\limits_{i=1}^k m_i^2 -1 \\[1ex]
& 1+r y_1 - (r+1)y_2 \geq 0  \\[1ex]  
& 1+s y_1-(s+1)y_2 \geq 0 \\[1ex] 
& y_1 \geq0, ~y_2 \geq 0.
\end{array}
\end{equation}
It is straightforward to derive a closed form expression for the optimal objective value of \eqref{LPsrg}  that is given in the following theorem.
\begin{thm} \label{strgClosedForm}
Let  $G=(V,E)$ be a  strongly regular graph with parameters $(n,\kappa,\lambda,\mu)$ and restricted eigenvalues $r \geq
0$ and $s<0$. Let $k$ and $m_i$ ($i=1,\ldots,k$) be positive integers such that $\sum_{j=1}^k m_j=n$. Then the SDP bound
${\rm GPP}_{\rm m}$ for the minimum GPP of $G$ is given by
\[
\max \left \{\frac{\kappa-r}{n} \sum_{i<j} m_i m_j, ~\frac{1}{2} \left ( n(\kappa+1) - \sum_i m_i^2 \right)  \right \}.
\]
Similarly, the  SDP bound ${\rm GPP}_{\rm m}$ for the maximum GPP is given by
\[
\min \left \{ \frac{\kappa -s}{n} \sum_{i<j} m_i m_j, ~\frac{1}{2}\kappa n \right \}.
\]
\end{thm}
We remark that in general, the point where this minimum ${\rm GPP}_{\rm m}$ is attained is the intersection point of the first constraint and the boundary of the third constraint, with objective value $\frac{\kappa-r}{n} \sum_{i<j} m_i m_j$. However, if $\frac{1}{n}\sum\limits_{i=1}^k m_i^2 -1 \leq -\kappa/s$, then the minimum is attained at the $y_1$-axis, with objective value $\frac{1}{2} ( n(\kappa+1) - \sum_i m_i^2)$.
For the case of the GEP, the results of Theorem \ref{strgClosedForm} coincide with the results from
\cite{deKlSotNaTr:10} and \cite{dKPaDoSo:10}. To see that, one should use the equation $n(\kappa +
rs)=(\kappa -s)(\kappa -r)$ (which follows from taking row sums of the equation $(A-rI)(A-sI)=(\kappa+rs) J$)
and other standard equations for the parameters of strongly regular graphs (see e.g., \cite{BrHa}), and the
fact that  ${\rm GPP}_{\rm m}$ is equivalent to the SDP relaxation for the GEP problem by Karish and Rendl
\cite{KarRend:98} (see \cite{Sot11}).

Next, we consider ${\rm GPP}_{\rm m\triangle}$ for SRGs.
From \eqref{CohTriangAgg} and \eqref{strY}, with  $A_1:=A$ and $A_2:=J-A-I$,
 it follows that for given $i,j,h \in \{1, 2\}$   the aggregated triangle inequality reduces to
\[
(\tr BA)y_1 + (\tr B(J-A-I))y_2 \leq b,
\]
where $B=p^i_{hj}A_i + p^h_{ij} A_h- p^j_{ih} A_{j}$ and  $b=p^i_{hj} \tr A_i J$.
After simplifying and removing equivalent inequalities, at most the inequalities
\begin{eqnarray}
y_1\leq 1, \quad y_2\leq 1,\label{tr1}\\
1+y_1 - 2y_2 \geq 0,\label{tr2}\\
1-2y_1 + y_2 \geq 0, \label{tr3}
\end{eqnarray}
remain, and when some of the intersection numbers $p^h_{ij}$ vanish, even fewer remain (we omit details for the sake of readability).
It is not hard to see that the constraints \eqref{tr1} are always redundant to the constraints of \eqref{LPsrg} (for example by drawing the feasible region),
and that \eqref{tr2} (cf.~\eqref{lineq2}) and \eqref{tr3} (cf.~\eqref{lineq3}) are redundant except for $r<1$ and $s>-2$, respectively,
which occurs only for the pentagon, disconnected SRGS, and complete multipartite graphs.
However, for the disconnected SRGs and complete multipartite graphs, the `nonredundant' constraints \eqref{tr2} and \eqref{tr3} (respectively)
don't occur precisely because of the vanishing of the relevant intersection numbers.
In other words, adding triangle inequalities to ${\rm GPP}_{\rm m}$ for strongly regular
graphs does {\em not} improve the bound, except possibly for the pentagon.
On the other hand, if we consider the pentagon and add the triangle inequalities to ${\rm GPP}_{\rm m}$ with $m=(2,3)^{\mathrm T}$, then the bound improves and is tight.

For the bisection problem  the aggregated independent set constraints are of the form \eqref{CohIndependAgg}, again with  $A_1:=A$ and $A_2:=J-A-I$.
Our numerical tests show that for many strongly regular graphs, the independent set constraints do not improve ${\rm GPP}_{\rm m}$,
but there are also graphs for which ${\rm GPP}_{\rm m-ind}$ dominates ${\rm GPP}_{\rm m}$, see Section \ref{sec:bisec} and \ref{sec:Aggregate}.

\subsection{Johnson and Kneser graphs} \label{sec:JohnKne}

In this section we show that for the Johnson and Kneser graphs (on triples), the SDP
bound  ${\rm GPP}_{\rm m}$ can be obtained by solving a linear programming problem. We also present
aggregated   triangle and independent set inequalities that one may add to  ${\rm GPP}_{\rm m}$. The Johnson
graphs were also studied by Karloff \cite{HK:99} in the context of the max cut problem,
in order to show that it is impossible to add valid linear inequalities to improve the performance ratio for
the celebrated Goemans-Williamson approximation algorithm. Our results show that ${\rm GPP}_{\rm m}$
improves after adding the independent set constraints.

The Johnson and Kneser graphs are defined as follows.
Let $\Omega$ be a fixed set of size $v$ and let $d$ be an integer such that $1\leq d\leq v/2$. The vertices of the Johnson scheme
$J(v,d)$ are the subsets of $\Omega$ with size $d$. The adjacency matrices of the association scheme are defined by the size of the intersection of these
subsets, in particular $(A_i)_{\omega,\omega'}=1$ if the subsets $\omega$ and $\omega'$ intersect in $d-i$ elements, for $i=0,\dots,d$. We remark that $A_1$
represents a so-called distance-regular graph $G$ --- the {\em Johnson graph} --- and $A_i$ represents being at distance $i$ in $G$.
The {\em Kneser graph} $K(v,d)$ is the graph with adjacency matrix $A_d$, that is, two subsets are adjacent whenever they are
disjoint. The Kneser graph $K(5,2)$ is the well-known Petersen graph.

For the case  $d=2$, the Johnson graph is strongly regular and also known as a triangular graph. Consequently
the  bound ${\rm GPP}_{\rm m}$ of $J(v,2)$  has a closed form expression (apply  Theorem \ref{strgClosedForm}
with $\kappa=2(v-2)$, $r=v-4$, and $s=-2$). Similarly,  the Kneser graph  $K(v,2)$  is strongly regular and
the closed form expression for the GPP follows from Theorem \ref{strgClosedForm} with
$\kappa=\binom{v}{2}-1-2(v-2)$, $r=1$, and $s=3-v$.

Here we focus on the next interesting group of Johnson and Kneser graphs, i.e., those on
triples ($d=3$), but we also note that the restriction to the case $d=3$ is not essential. The eigenvalues
(character table) of the Johnson scheme can be expressed in terms of Eberlein polynomials; see Delsarte's
thesis \cite[Thm.~4.6]{Delsarte73}. For $d=3$, the character table is
\begin{equation} \label{charTabJoh}
P= \left (
\begin{array}{cccl}
1 & \theta_0 & \varphi(\theta_0) & \binom{v}{3}-1 -\theta_0-\varphi(\theta_0) \\[1ex]
1 & \theta_1 & \varphi(\theta_1) & -1 -\theta_1-\varphi(\theta_1)\\[1ex]
1 & \theta_2 & \varphi(\theta_2) & -1 -\theta_2-\varphi(\theta_2)  \\[1ex]
1 & \theta_3 & \varphi(\theta_3) & -1 -\theta_3-\varphi(\theta_3)
\end{array}
\right ),
\end{equation}
where $\varphi(\theta) = \frac{1}{4} \left (\theta^2-(v-2)\theta -3(v-3) \right) $, $\theta_0 = 3(v-3)$,
$\theta_1 = 2v-9$, $\theta_2=v-7$, and $\theta_3 = -3$.

Let $A_1$ denote the adjacency matrix of $J(v,3)$,
and $A_3$ the adjacency matrix of $K(v,3)$. We first simplify  ${\rm GPP}_{\rm m}$  for the case that the
graph under consideration is the Johnson graph $J(v,3)$. From Theorem \ref{thmAlg}, it follows that there
exists an optimal solution $Y$ to ${\rm GPP}_{\rm m}$ which belongs to the coherent algebra spanned by $\{I,
A_1, A_2, A_3 \}$.  Thus, there exist $y_1,y_2, y_3 \geq 0$ such that
\begin{equation}\label{strYJohn}
Y= I + y_1A_1 + y_2A_2 + y_3A_3.
\end{equation}
Now, similar to the case of strongly regular graphs (see also \cite{dKPaDoSo:10,deKlSotNaTr:10}), we can
rewrite the objective function and constraints from ${\rm GPP}_{\rm m}$  by using \eqref{strYJohn}. The
derived LP for the minimum GPP is
\begin{equation} \label{LPJohn}
\begin{array}{ll}
\min &   \frac{3}{2}\binom{v}{3}(v-3) (1-y_1) \\[1ex]
{\rm s.t.} &  A_{eq}y = b_{eq} \\[1ex]
&   A_{neq} y \geq b_{neq} \\[1ex]
&  y\geq 0, ~y\in \R^3,
\end{array}
\end{equation}
where
\begin{equation}\label{Aeq}
A_{eq} := \left ( 3(v-3), ~3\binom{v-3}{2}, ~\binom{v-3}{3} \right ),
\end{equation}

\begin{equation} \label{Aineq}
A_{neq}
:= \left (
\begin{array}{ccc}
2v-9 & \frac{1}{2} (v^2-13v +36) & \frac{1}{2} (-v^2+9v-20) \\[1.5ex]
v-7 & -2v+11 & v-5  \\[1.5ex]
-3 & 3 & -1
\end{array} \right ),
\end{equation}
\begin{equation} \label{beq}
b_{eq}:=\frac{1}{n} \sum\limits_{i=1}^k m_i^2 -1, \quad
b_{neq} := -(1, 1, 1)^{\mathrm{T}},
\end{equation}
and $n=\binom{v}{3}$ is the number of vertices of $J(v,3)$. To derive (\ref{LPJohn}) we exploited the fact
that the matrices $\{I, A_1, A_2, A_3 \}$  may be simultaneously diagonalized and we used the character table
\eqref{charTabJoh}. Note that the computation time for solving \eqref{LPJohn} is negligible and does not
increase with the order of the Johnson graph.

\begin{thm} \label{JohnsonLP}
Let  $J(v,3)$ be the Johnson graph, with $n$ vertices, and let $k$ and $m_i$ ($i=1,\ldots,k$) be positive
integers such that $\sum_{i=1}^k m_i=n$. Then the SDP  bound ${\rm GPP}_{\rm m}$  for the minimum GPP of
$J(v,3)$ is equal to the optimal value of the linear programming problem \eqref{LPJohn}.
\end{thm}

Similarly, we  simplify ${\rm GPP}_{\rm m}$  for the Kneser graph $K(v,3)$. Clearly, the only difference is
the objective function which corresponds to the partition of the Kneser graph. The resulting LP relaxation is
\begin{equation} \label{LPKnes}
\begin{array}{ll}
\min &   \frac{1}{2}\binom{v}{3}\binom{v-3}{3} (1-y_3) \\[1ex]
{\rm s.t.} &  A_{eq}y = b_{eq} \\[1ex]
&   A_{neq} y \geq b_{neq} \\[1ex]
&  y\geq 0, ~y\in \R^3,
\end{array}
\end{equation}
where $A_{eq}$-$b_{neq}$ are as in \eqref{Aeq}-\eqref{beq}. This leads to the following result.

\begin{thm} \label{KneserLP}
Let  $K(v,3)$ be the Kneser graph on $n$ vertices, and let $k$ and $m_i$ ($i=1,\ldots,k$) be
 positive  integers such that $\sum_{i=1}^k m_i=n$.
Then the SDP  bound ${\rm GPP}_{\rm m}$  for the minimum GPP of $K(v,3)$  is equal to the optimal value of the
linear programming problem \eqref{LPKnes}.
\end{thm}

We can add to \eqref{LPJohn} and \eqref{LPKnes} the aggregated triangle inequalities \eqref{CohTriangAgg}.
For given $i,j,h \in \{1, 2, 3\}$, these reduce to
\[
(\tr BA_1)y_1 + (\tr BA_2)y_2 + (\tr BA_3)y_3 \leq b,
\]
where  $B=p^i_{hj}A_i + p^h_{ij} A_h- p^j_{ih} A_{j}$ and  $b=p^i_{hj} \tr A_i J$. After
taking into consideration all possible choices of $i,j,h \in \{1, 2, 3\}$, there remain only seven
(aggregated) triangle inequalities when $v=6$ and eleven when $v>6$. Numerical results
indicate that these additional inequalities do not improve the solution obtained by
solving \eqref{LPJohn} and \eqref{LPKnes}.

Since we know how to aggregate the independent set constraints \eqref{CohIndependAgg} when $k=2$,
we tested the effect on the bound ${\rm GPP}_{\rm m}$ of adding these constraints. The numerical results show that the bound may improve, e.g.,
 for the bisection of $J(7,3)$ with $m=(17,18)^{\mathrm T}$, it improves from $62$ to $64$.

\section{A new eigenvalue bound for the GPP} \label{sec:AnyGaph}

In this section we present a closed form expression for the optimal value of a relaxation for the GPP for any graph and any $k\geq 2$.
To the best of our knowledge, in the literature there are such general closed form bounds for the GPP  only for
 $k=2$ (see e.g., Juvan and Mohar \cite{JuMo:92} and Falkner, Rendl, and Wolkowicz \cite{FeReWo:92}) and $k=3$ (see \cite{FeReWo:92}).
 Recently, Pong at al.~\cite{PoSuWaWolkowicz} and Rendl et al. \cite{RendlLissPia} derived eigenvalue bounds for the vertex separator problem, which is closely
 related to the GPP.

In order to derive an eigenvalue bound for the GPP, we  relax several constraints in ${\rm GPP}_{\rm m}$.
In particular, we relax $\diag(Y)=u_n$ to $\tr Y=n$ and remove nonnegativity constraints.
Moreover, we use \eqref{obj1} to rewrite the objective in terms of the Laplacian matrix $L$, which leads to the relaxation
\begin{equation} \label{RSrelax}
\begin{array}{rl}
\min & \frac{1}{2} \tr LY \\[1ex]
{\rm s.t.} & \tr Y=n \\
& \tr JY=\sum\limits_{i=1}^k m_i^2  \\[1.5ex]
& kY - J_n\succeq 0.
\end{array}
\end{equation}
Recall from Section \ref{Ex:LaplAlg} that we denote the distinct Laplacian eigenvalues of the graph by $0=
\lambda_0 < \lambda_1 < \ldots < \lambda_d$, and their corresponding multiplicities
$f_i$, for $i=0,\ldots,d$, and let $\mathcal L=\spen \{ F_0,\ldots,F_d\}$ be the Laplacian algebra of the
graph. By Theorem \ref{thmAlg}, there exists an optimal solution $Y$ to \eqref{RSrelax} in $\mathcal L$, and
therefore we may assume that
$Y=\sum_{i=0}^d y_iF_i,$
where $y_i\in  \R$ ($i=0,\ldots,d$) (as before, these are the new
variables). We will exploit this to rewrite \eqref{RSrelax}. The objective is
\[
\tr LY= \tr  (\sum_{i=0}^d \lambda_i F_i )( \sum_{j=0}^d y_j F_j )    =\sum_{i=0}^d \lambda_if_iy_i.
\]
The constraint $ \tr Y=n $ reduces to $ \sum_{i=0}^d f_i y_i = n, $ while the constraint $\tr JY=\sum_{i=1}^k
m_i^2 $ reduces to $y_0=  (\sum_{i=1}^k m_i^2)/n$. It remains only to reformulate the semidefinite
constraint:
\[ kY - J =k\sum_{i=0}^d y_iF_i - J=
\frac{k\sum_{i=1}^k m_i^2  -n^2 }{n^2}J + \sum_{i=1}^d y_i F_i \succeq 0.
\]
From this, it follows that $y_i\geq 0$. To conclude, the SDP relaxation \eqref{RSrelax}
reduces to
\begin{equation}\label{LPGraph}
\begin{array}{rl}
\min & \frac12 \sum\limits_{i=1}^d \lambda_if_iy_i \\[2ex]
{\rm s.t.} & \sum\limits_{i=1}^d f_i y_i = \frac{2}{n} \sum\limits_{i<j}m_im_j  \\[2ex]
& y_i \geq 0, \quad  i=1,\ldots,d.
\end{array}
\end{equation}

\begin{thm} \label{thm:AnyGraph}
Let  $G$ be a graph on $n$ vertices, and let  $k$ and $m_i$ $(i=1,\ldots,k)$ be positive integers such that
$\sum_{i=1}^k m_i=n$. Then the SDP lower bound \eqref{RSrelax} for the minimum GPP of $G$ is equal to
\[
\frac{\lambda_1}{n} \sum\limits_{i<j}m_im_j,
\]
and the SDP upper  bound for the maximum GPP of $G$  that is obtained by replacing
$\min$ by $\max$ in \eqref{RSrelax} is
\[
\frac{\lambda_d}{n} \sum\limits_{i<j}m_im_j.
\]
\end{thm}

\noindent
{\em Proof}. This follows from \eqref{LPGraph}. \qed \\

\noindent Our  eigenvalue bound for the bisection problem (the case $k=2$)
coincides with a well-known result in spectral graph theory, see \cite{JuMo:92,MohPolj:93}.
 Therefore, Theorem \ref{thm:AnyGraph} may be seen as a generalization of this
result for the $2$-partition problem to any $k$-partition problem.  Falkner, Rendl, and Wolkowicz
\cite{FeReWo:92} derived a closed form bound for the minimum $3$-partition problem of the form
\begin{equation}\label{FRW}
\frac{1}{2} \theta_1 \mu_1 + \frac{1}{2} \theta_2 \mu_2,
\end{equation}
where $\mu_{1,2}= (m_1m_2+m_1m_3+m_2m_3 \pm \sqrt{ m_1^2m_2^2 + m_1^2m_3^2  +  m_2^2m_3^2  -nm_1m_2m_3 })/n$,
and $\theta_1$ and $\theta_2$ are the two smallest nonzero (not necessarily distinct)
Laplacian eigenvalues. It is clear that this lower bound coincides with ours when $\theta_1=\theta_2$
($=\lambda_1$). To the best of our knowledge there are no other closed form bounds
for the minimum $k$-partition problem when $k>3$, or for the maximum $k$-partition problem when $k> 2$
that is applicable to any graph.
Although the bounds from Theorem \ref{LPGraph} are, in general, dominated by the bounds obtained from ${\rm GPP}_{\rm m}$,
they may be useful in the theoretical analysis of the GPP and related problems.
 Still, our numerical results show that for many problems
the new eigenvalue bound is equal to ${\rm GPP}_{\rm m}$, see Section \ref{sec:NumRez}. We finally remark
that for strongly regular graphs, the eigenvalue bounds also follow from Theorem \ref{strgClosedForm} because
$\lambda_1=\kappa-r$ and $\lambda_d=\kappa-s$ (indeed, they are even the same unless $\sum_i m_i^2 <
\frac{r+1}{n-\kappa+r}n^2$ or $\sum_i m_i^2 < \frac{-s}{\kappa-s}n^2$, respectively). This is related to the
fact that the Laplacian algebra and the used coherent algebra are the same for strongly regular graphs.

\section{Improved relaxations for the GPP} \label{sect:EQP}

\subsection{An improved relaxation from the quadratic assignment problem} \label{sec:impVec}

In this section, we derive a new  SDP  relaxation for the GPP that is obtained by strengthening  the SDP
relaxation of the more general quadratic assignment problem by  Zhao, Karisch, Rendl, and Wolkowicz
\cite{ZhKaReWo:98}, by adding two constraints that correspond to assigning two vertices of the graph to
different parts of the partition. A similar approach was used by the authors \cite{vDamSo:12} to derive the
best known bounds for the bandwidth problem of graphs with symmetry. For the equipartition problem,
the new relaxation dominates the relaxation for the GEP by De Klerk et al.~\cite{dKPaDoSo:10}, that
is obtained by fixing one vertex of the graph. Our new bound is however not restricted to the equipartition
problem, and it is also suitable for graphs with symmetry.

The GPP is a special case of the quadratic assignment problem
\[
\min_{X\in \Pi_n} \frac{1}{2} \tr AXBX^{\mathrm T},
\]
where $A$ and $B$ are given symmetric $n \times n$ matrices, and $\Pi_n$  is the set of $n\times n$ permutation matrices.
For the graph partition problem, $A$ is the adjacency matrix of the relevant graph $G$ with $n$ vertices, and $B$ is the adjacency matrix of
the complete multipartite graph $K_{m_1,\ldots,m_k}$ with  $k$ classes of sizes $m_1$,\ldots,$m_k$ (with $m_1+\ldots +m_k=n$).
For example, for  the $k$-equipartition problem with $n=km$,
\begin{equation}\label{defBEQ}
B:=(J_k-I_k)\otimes J_m,
\end{equation}
(where $\otimes$ is the Kronecker product). In the general case, $B$ has the same block structure, but the
sizes $m_1$,\ldots,$m_k$ of the blocks vary. In particular, for the bisection problem with
$m=(m_1,m_2)^{\mathrm T}$,
\begin{equation} \label{defBBis}
B:=\left (
\begin{array}{cc}
0_{m_1\times m_1} & J_{m_1\times m_2}\\
J_{m_2\times m_1} & 0_{m_2\times m_2}\\
\end{array} \right ).
\end{equation}
\noindent Now it follows that the following `vector-lifting' SDP relaxation of this particular
QAP (see Zhao et al.~\cite{ZhKaReWo:98} and Povh and Rendl \cite{PoRe:09}) is  also a relaxation
for the GPP:
\[
(\mbox{GPP}_{\rm QAP})~~~~
\begin{array}{rcl}
\min && \frac{1}{2} \tr(B\otimes A)Y\\[1ex]
{\rm s.t.} && \tr(I_n\otimes E_{jj})Y=1, ~~\tr(E_{jj}\otimes I_n)Y=1, \quad  j=1,\ldots,n\\[1ex]
&& \tr (I_n\otimes(J_n-I_n)+(J_n-I_n)\otimes I_n)Y=0\\[1ex]
&& \tr JY=n^2 \\[1ex]
&& Y\geq 0, ~~Y \succeq 0,
\end{array}
\]
where (here and below) $E_{ij}=e_i e_j^{\mathrm{T}}$. In \cite{Sot10} (see also \cite{Sot11}) it is proven
that for the equipartition  problem, the relaxations $\mbox{GPP}_{\rm QAP}$ and $\mbox{GPP}_{\rm m}$ are
equivalent, and in \cite{Sot11} that the first dominates the second for the bisection problem. De Klerk et
al.~\cite{dKPaDoSo:10} strengthened  $\mbox{GPP}_{\rm QAP}$ for the GEP  by adding a
constraint that corresponds to assigning an arbitrary vertex of the complete multipartite graph to a vertex
in the graph.

Here, we extend the approach from \cite{dKPaDoSo:10} (see also  \cite{vDamSo:12}) and assign (several times) a pair of vertices of $G$ to an edge in
$K_{m_1,\ldots,m_k}$. By symmetry, we have to do this for one pair of vertices in each orbital (recall from Section \ref{sec:SyminGr}
that the orbitals  actually  represent the `different' kinds of pairs of vertices; (ordered) edges, and (ordered) nonedges  in the graph $G$).
Let us assume that there are $t$ such orbitals ${\mathcal O}_h$ ($h=1,2\ldots,t$) of edges and nonedges, and
note that for highly symmetric graphs, $t$ is relatively small.
We formally state the above idea in the following theorem.
\begin{thm} \label{thmFixnew}
Let $G$ be an undirected graph on $n$ vertices with adjacency matrix $A$, and let $\mathcal{O}_h$
$(h=1,2,\dots,t)$ be the orbitals of edges and nonedges coming from the automorphism group of $G$. Let
$(s_1,s_2)$ be an arbitrary edge in  $K_{m_1,\ldots,m_k}$ while $(r_{h1},r_{h2})$ is an arbitrary pair of
vertices in $\mathcal{O}_h$ $(h=1,2,\dots,t)$. Let $\Pi_n(h)$ be the set of matrices $X \in \Pi_n$ such that
$X_{r_{h1},s_1}=1$ and $X_{r_{h2},s_2}=1$ $(h=1,2,\dots,t)$. Then
\[
\min_{X \in \Pi_n} \tr X^{\mathrm{T}}AXB =
\min_{h=1,\dots,t} \min_{X \in \Pi_{n}(h)} \tr X^{\mathrm{T}}AXB.
\]
\end{thm}

\proof Similar to  the proof of Theorem 10 in \cite{vDamSo:12}. \qed \\

Clearly, exploiting this requires solving several SDP (sub)problems. However, if we assign to an edge
$(s_1,s_2)$  in  $K_{m_1,\ldots,m_k}$ a pair of vertices  $(r_{h1},r_{h2})$ from  $\mathcal{O}_h$
$(h=1,2,\dots,t)$, then we can add to $\mbox{GPP}_{\rm QAP}$ the constraints
\[
\tr (E_{s_i,s_i} \otimes E_{r_{hi},r_{hi}})Y=1, \quad i=1,2.
\]
Thus, we obtain several SDP problems of the form
\begin{equation} \label{ZWGEPh}
\begin{array}{rcl}
\mu^h := \min && \frac{1}{2} \tr(B\otimes A)Y\\[1ex]
{\rm s.t.} && \tr(I_n\otimes E_{jj})Y=1, ~~\tr(E_{jj}\otimes I_n)Y=1, \quad  j=1,\ldots,n\\[1ex]
&& \tr (E_{s_i,s_i} \otimes E_{r_{hi},r_{hi}})Y=1, \quad i=1,2 \\[1ex]
&& \tr (I_n\otimes(J_n-I_n)+(J_n-I_n)\otimes I_n)Y=0\\[1ex]
&& \tr JY=n^2 \\[1ex]
&& Y\geq 0, ~~Y \succeq 0,
\end{array}
\end{equation}
where $h=1,\ldots,t$, and the new lower bound for the GPP is \label{relaxFix}
\[
{\rm GPP}_{\rm fix}=\min\limits_{h=1,\ldots,t}\mu^h.
\]
We remark that $\mu^h$ is a relaxation that depends on the particular edge $(s_1,s_2)$ (but not on the
particular pair $(r_{h1},r_{h2}) \in \mathcal{O}_h$). However, for the GEP and GBP (but not in general!), the
lower bound ${\rm GPP}_{\rm fix}$ is independent of the edge $(s_1,s_2)$. This is due to the fact
that $K_{m,\ldots,m}$ and $K_{m_1,m_2}$ are edge-transitive.

The following proposition follows directly from  \eqref{ZWGEPh}.
\begin{cor}\label{cor:newbounddominates}
Let $(s_1,s_2)$ be an arbitrary edge in  $K_{m_1,\ldots,m_k}$ and  $(r_{h1},r_{h2})$ be an arbitrary pair of
vertices in $\mathcal{O}_h$ $(h=1,2,\dots,t)$. Then the SDP bound  ${\rm GPP}_{\rm fix}$ dominates
${\rm GPP}_{\rm QAP}$.
\end{cor}
Similarly, the following corollary for the GEP follows.
\begin{cor}\label{cor:newbounddominates1}
Consider the equipartition problem. Let $(s_1,s_2)$ be an arbitrary edge in
$K_{m,\ldots,m}$, and $(r_{h1},r_{h2})$ be an arbitrary pair of vertices in $\mathcal{O}_h$
$(h=1,2,\dots,t)$. Then the SDP relaxation ${\rm GPP}_{\rm fix}$ dominates the SDP relaxation from
\cite[Eq. 10]{dKPaDoSo:10}.
\end{cor}

It is, in general, hard to solve \eqref{ZWGEPh} (and thus  ${\rm GPP}_{\rm fix}$) for $n\geq 16$, see e.g., \cite{ReSo}.
Therefore we need to further exploit the symmetry of $K_{m_1,\ldots,m_k}$ (in particular, consider pointwise stabilizers)
and the graphs under consideration, see also \cite{vDamSo:12}.
We do this by applying the general theory of symmetry reduction to the SDP subproblems \eqref{ZWGEPh}  in a mechanical way, as described in,
e.g., \cite{vDamSo:12,deKlSot:10,dKPaDoSo:10,dKSo:12}.
After symmetry reduction of \eqref{ZWGEPh} the largest linear matrix inequality contains matrices of size $3n$ (resp.~$2n$) for the GEP (resp.~GBP).

Our numerical results show that ${\rm GPP}_{\rm fix}$ can be a significantly stronger bound than ${\rm GPP}_{\rm m}$ for highly
symmetric graphs (i.e., for which $t$ is very small) and for cases that the bound obtained by solving
${\rm GPP}_{\rm m}$ cannot be improved by adding triangle and independent set constraints. This could be a
consequence of the fact that (some) symmetry in the graph has been broken.

\subsection{An improved matrix-lifting relaxation} \label{sec:impMat}

Clearly, we can exploit the idea of fixing a pair of vertices in a graph also in the context of the  matrix lifting relaxation  ${\rm GPP}_{\rm m}$.
Assume again that for the given graph $G$ there are $t$ orbitals ${\mathcal O}_h$ ($h=1,2\ldots,t$) of edges and nonedges.
Now, in order to assign two (arbitrary) vertices $(r_{h1},r_{h2})\in {\mathcal O}_h$  of the graph $G$
to two different subsets, we add to ${\rm GPP}_{\rm m}$ the constraint
\[
\tr ( E_{r_{h1},r_{h2}}+E_{r_{h2},r_{h1}})Y=0.
\]
Therefore, computing this new lower bound reduces to solving $t$ subproblems of the form
\begin{equation} \label{RSfixPom}
\begin{array}{rl}
\nu^*_h :=\min & \frac{1}{2} \tr A(J_n -Y) \\[1ex]
{\rm s.t.} & \diag(Y)=u_n \\[1ex]
& \tr JY=\sum\limits_{i=1}^k m_i^2  \\[2ex]
& \tr ( E_{r_{h1},r_{h2}}+E_{r_{h2},r_{h1}})Y=0 \\[2ex]
& kY - J_n\succeq 0, ~~Y\geq 0,
\end{array}
\end{equation}
($h=1,\ldots,t$).
Consequently, the new matrix-lifting lower bound is a minimum over $t$ SDP bounds, i.e.,
\begin{equation} \label{RSfix}
\min\limits_{h=1,\ldots, t} \nu^*_h.
\end{equation}
The following result follows immediately.
\begin{cor}
The SDP bound \eqref{RSfix} dominates ${\rm GPP}_{\rm m}$.
\end{cor}
\noindent
In order to solve \eqref{RSfixPom} (and thus \eqref{RSfix}) we further exploit symmetry in the graphs under consideration
in a similar way as described in Section \ref{sec:aggregateConstr}.

Our numerical results suggest that the new SDP  bound \eqref{RSfix} is dominated by
$\mbox{GPP}_{\rm fix}$, and also that \eqref{RSfix} is not dominated by $\mbox{GPP}_{\rm QAP}$, or vice versa.

\section{Numerical results} \label{sec:NumRez}

In this section we present numerical results for the graph partition problem.
In particular, we compare bounds from all the presented relaxations and several relaxations from the literature.
All relaxations were solved with SeDuMi \cite{sedumi} using the Yalmip interface \cite{YALMIP} on
an Intel Xeon X5680, $3.33$ GHz dual-core processor with 32 GB memory. To compute orbitals, we used GAP \cite{gap}.

\subsection{Why symmetry?}

We first show the importance of exploiting symmetry in graphs, when applicable, in order to compute SDP
bounds. In Table \ref{tab1} we consider the planar unweighted {\tt grid graphs}, where $|V|= \sharp$ rows
$\times$ $\sharp$ columns. They are generated by the rudy graph generator \cite{Rinaldi}. Table \ref{tab1}
presents computational times, in seconds, required to solve ${\rm GPP}_{\rm m}$ with and without exploiting
symmetry (see also Table $3$ in \cite[online supplement]{Sot11}). The table reads as
follows. In the first two columns, the sizes of the graphs and the sizes of the partitions are
specified. The third column lists computational times required to solve ${\rm GPP}_{\rm m}$ without
exploiting symmetry. The fourth column provides the rank of the associated matrix $*$-algebra that is
obtained as the centralizer ring of the automorphism group of the graph, and the last column contains
computational times required to solve ${\rm GPP}_{\rm m}$ after exploiting symmetry. Note that even though
the graphs are not highly symmetric, the reduction in computational times after exploiting symmetry is
significant.

\begin{table}[h!]
\caption{Computational time (s.) to solve ${\rm GPP}_{\rm m}$ for the min 3-partition problem.}  \label{tab1}
\begin{center}
\begin{tabular}{ccccc}
\hline\noalign{\smallskip}
$|V|$ & $m^{\mathrm T}$ & no symmetry  &  $r_{\rm aut}$ & symmetry\\[1ex]
\noalign{\smallskip}\hline\noalign{\smallskip}
$9 \times 9$ & $(35, 30, 16)$ & 198.35 & 861 & 1.70 \\
$10 \times 10$ & $(50, 25, 25)$ & 799.21 & 1275 & 3.41 \\
\noalign{\smallskip}\hline
\end{tabular}
\end{center}
\end{table}

\subsection{Combinatorial symmetry vs. group symmetry}\label{ces:CombSymNum}

In this section we list numerical results for several graphs that have (substantially) more
combinatorial symmetry than symmetry coming from the automorphism group. We provide the eigenvalue bound of
Theorem \ref{thm:AnyGraph} and ${\rm GPP}_{\rm m}$ for the GEP of those graphs.

Table \ref{tab:other} reads as follows. In the first three columns, we list the graphs, the number of
vertices, and the number of parts $k$ of the equipartition, respectively. {\tt Chang3} is one of the
strongly regular graphs introduced by Chang \cite{Chang} (see also \cite{BrChang}). For a
description of the {\tt Doob} graph, see Section \ref{sec:comb}. Graphs {\tt A64v30} and {\tt A64vEnd}  are
strongly regular graphs with parameters $(64,18,2,6)$ obtained by Haemers and Spence \cite{HaeSpe:01}, where
30 (resp.~End) means that it is the 30th (resp.~last) graph in the list, see also
\url{http://www.maths.gla.ac.uk/~es/SRGs/64-18-2-6}. The {\tt design} graph is the bipartite incidence graph
of a symmetric $2$-$(45,12,3)$-design, see \url{http://www.maths.gla.ac.uk/~es/polar/45-12-3.36}. In the
fourth column of Table \ref{tab:other} we give the eigenvalue lower bound from Theorem \ref{thm:AnyGraph},
while in  the fifth column, we list the SDP bound ${\rm GPP}_{\rm m}$. All presented bounds are rounded up to
the closest integer. In the last four columns we list the rank of the  coherent configuration
corresponding to the graph's combinatorial symmetry, computational times required to solve ${\rm GPP}_{\rm
m}$ after exploiting its combinatorial symmetry,
 the rank of the  coherent configuration coming from the automorphism group, and the corresponding computational times,
 respectively.
If a graph is strongly regular we do not report the computational time since for such graphs we use the closed form expression from  Theorem \ref{strgClosedForm}.

Note that  {\tt A64vEnd} does not have any symmetry coming from automorphisms, but it has lots of combinatorial
symmetry. We remark that the listed graphs are {\em not} isolated cases, but only a sample that shows that
combinatorial symmetry may differ significantly from group symmetry.

\begin{table}[h!]
\caption{Lower bounds and computational times (s.) for the min GEP. }  \label{tab:other}
\begin{center}
\begin{tabular}{ccccccccr}
\hline\noalign{\smallskip}
$G$ & $n$ & $k$  & eig & ${\rm GPP}_{\rm m}$ & $r_{\rm comb}$ & time  & $r_{\rm aut}$ & time \\
\noalign{\smallskip}\hline\noalign{\smallskip}
{\tt Chang3} & 28 & 7 & 96 & 126 & 3 & -- &  14 & 0.23 \\
{\tt A64v30} & 64 & $8$ & 448 & 448 & 3 & -- & 90 & 0.61 \\
{\tt Doob} & 64 & $8$ & 112 & 160 & 4 &  0.34 & 8 & 0.41 \\
{\tt A64vEnd} & 64 & $4$ & 384 & 384 & 3 &  -- & -- & 14.33  \\
{\tt design} & 90 & $9$ & 360 & 360 & 4 & 0.40 & 2074 & 4.56 \\
\noalign{\smallskip}\hline
\end{tabular}
\end{center}
\end{table}

\subsection{The graph equipartition problem}

In this section we compare different relaxations for the  equipartition problem. We first present results for
the {\tt Higman-Sims} graph \cite{HigSim:68}, see Table \ref{tab2}. The {\tt Higman-Sims} graph  is a
strongly regular graph  with parameters $(100,22,0,6)$. The max and min $k$-equipartition problem
for this graph was studied in \cite{dKPaDoSo:10,deKlSotNaTr:10}.

Table \ref{tab2} reads as follows. The first
column specifies whether we are solving a minimization or maximization problem, while the second
column shows the number of parts $k$ of the equipartition. The third column provides the new eigenvalue
bound, see Theorem \ref{thm:AnyGraph}. The fourth column lists $\mbox{GPP}_{\rm m}$ which is known
to be equivalent to  $\mbox{GPP}_{\rm QAP}$ for the case of the equipartition (for a proof, see
\cite{Sot10,Sot11}). The fifth column provides bounds obtained by solving the relaxation from
\cite{dKPaDoSo:10} (that is, the improved  $\mbox{GPP}_{\rm QAP}$  by adding  a constraint that
corresponds to fixing a vertex in the graph). In the sixth column, we list the bounds obtained by solving
GPP$_{\rm fix}$, see page \pageref{relaxFix}, while the seventh column contains the bounds obtained by
solving the SDP relaxation from  \cite{deKlSotNaTr:10}. The latter relaxation is the `level two
reformulation-linearization technique-type relaxation for the QAP with an additional linear matrix inequality
constraint', which is known to be at least as strong as $\mbox{GPP}_{\rm QAP}$.

In  Table \ref{tab2}, the bounds   improve along with
 increasing complexity of the relaxations; the strongest bound is from \cite{deKlSotNaTr:10}. Note however
that the bound from \cite{deKlSotNaTr:10} is appropriate only for vertex-transitive graphs, while our bounds
do {\em not} have such a restriction. The last column provides bounds obtained from heuristics that are taken from
Table $2$ and $3$ in \cite{deKlSotNaTr:10}.

We remark that the SDP bound \eqref{RSfix} provides the same bounds as  $\mbox{GPP}_{\rm m}$ for all problems in Table \ref{tab2}.\\

\begin{table}[h!] 
\caption{Bounds for the GEP for the {\tt Higman-Sims} graph.}  \label{tab2}
\begin{center}
\begin{tabular}{lccccccc}
\hline\noalign{\smallskip}
$\max$ & $k$ & eig & $\mbox{GPP}_{\rm m}$ & \cite{dKPaDoSo:10}  &  GPP$_{\rm fix}$ & \cite{deKlSotNaTr:10} & lower bound\\
\noalign{\smallskip}\hline\noalign{\smallskip}
 & 4  & 1125 & 1100 & 1097 & 1094 & 1048 & 1006\\
 & 5 & 1200 & 1100 & 1100 & 1100& 1100& 1068\\
 \noalign{\smallskip}\hline\hline\noalign{\smallskip}
$\min $ & &&&&& & upper bound\\
\noalign{\smallskip}\hline\noalign{\smallskip}
& 20 & 950 & 950 & 951 & 951 & 975 & 980 \\
& 25 & 960 & 960 & 963 & 964 & 1000 & 1000 \\
\noalign{\smallskip}\hline
\end{tabular}
\end{center}
\end{table}

In Table \ref{tab:maxKpart1} we present results for the maximum equipartition problem
for several Johnson and Kneser graphs (see Section \ref{sec:JohnKne}). The table reads as follows. In the first column we list the graphs,
in the second column the number of vertices, and in the third column the number of parts $k$ of the equipartition.
In the fourth and fifth column we list the new eigenvalue bound and ${\rm GPP}_{\rm m}$, respectively.
In the last two columns of Table \ref{tab:maxKpart1} we list  GPP$_{\rm fix}$ and the corresponding required computational times.

In Table \ref{tab:maxKpart1} we do not report  computational times for   the new eigenvalue bound and ${\rm GPP}_{\rm m}$ since they are negligible,
see Sections \ref{sec:SRG}, \ref{sec:JohnKne}, and \ref{sec:AnyGaph}.
Also, since adding triangle inequalities to  ${\rm GPP}_{\rm m}$ for the problems  in Table \ref{tab:maxKpart1} do not improve ${\rm GPP}_{\rm m}$,
we did not make a separate column for ${\rm GPP}_{\rm m\triangle}$.
All presented bounds are rounded down to the closest integer.

The results show that the new bound  GPP$_{\rm fix}$ can be significantly stronger than ${\rm GPP}_{\rm m}$,
in particular for problems when the eigenvalue bound and ${\rm GPP}_{\rm m}$ provide the same bound.
The results also show that the eigenvalue bound performs well for most of the instances.
\begin{table}[h!]
\caption{Upper bounds and computational times (s.)  for the max GEP.}  \label{tab:maxKpart1}
\begin{center}
\begin{tabular}{ccccccr}
\hline\noalign{\smallskip}
$G$ &$n$ & $k$ & eig & ${\rm GPP}_{\rm m}$ &  GPP$_{\rm fix}$ & time  \\
\noalign{\smallskip}\hline\noalign{\smallskip}
$K(8,2)$ & 28 & 4 & 210  &  210 & 204 & 7.66 \\
$K(9,2)$ & 36 & 3 & 324 & 324 &  317 & 14.43  \\
$K(9,2)$ & 36 & 12 &  444 &  378 & 378 & 5.32 \\
$K(12,2)$ & 66 & 6 & 1485 &  1485 & 1473 & 31.79 \\
$J(8,3)$ & 56 & 4  & 378  & 378 & 377 & 148.54 \\
$K(9,3)$ & 84 & 3 & 840 &  840   & 828 &  551.99 \\
$K(15,2)$ & 105 & 5 & 3780 & 3780 & 3772 &  106.97 \\
$K(10,3)$ & 120 & 3 & 2000 & 2000 & 1979 & 1097.10  \\
\noalign{\smallskip}\hline
\end{tabular}
\end{center}
\end{table}

\subsection{The graph bisection problem} \label{sec:bisec}
In this section we present numerical results for the graph bisection problem. All graphs in Table
\ref{tab:GBP} are strongly regular.
The table reads as follows. In the first column we list the graphs. The
{\tt Johnson} graphs are defined in Section \ref{sec:JohnKne}, whereas the {\tt Hoffman-Singleton} ({\tt HS})
graph, the {\tt Gewirtz} graph, and the $M_{22}$ graph are the unique strongly regular graphs with parameters
$(50,7,0,1)$, $(56,10,0,2)$, and $(77,16,0,4)$, respectively. In the second column of Table \ref{tab:GBP} we
list the number of vertices in the corresponding graph, while the third column contains the sizes of the
subsets. We choose these sizes arbitrarily. In the remaining columns we provide the lower bounds
$\mbox{GPP}_{\rm m}$, \eqref{RSfix}, $\mbox{GPP}_{\rm QAP}$,   $\mbox{GPP}_{\rm m-ind}$, and GPP$_{\rm fix}$,
respectively. All bounds are rounded up to the closest integer. In Table \ref{tab:GBPtime} we provide
computational times required to solve the problems from Table \ref{tab:GBP} (the times to compute
GPP$_{\rm fix}$ and \eqref{RSfix} are  sums of computational times of all subproblems needed to obtain the
bounds).

From Table \ref{tab:GBP}  it follows that \eqref{RSfix} is not dominated by $\mbox{GPP}_{\rm QAP}$, or vice
versa. Similarly, we may conclude that $\mbox{GPP}_{\rm QAP}$ and GPP$_{\rm fix}$ are not dominated by $\mbox{GPP}_{\rm m-ind}$,
or vice versa.
One more  interesting observation is that for all tested instances, the new eigenvalue bound is equal to the
bound obtained by solving $\mbox{GPP}_{\rm m}$.

\begin{table}[h!]
\caption{Lower bounds for the min GBP.}  \label{tab:GBP}
\begin{center}
\begin{tabular}{cccccccc}
\hline\noalign{\smallskip}
$G$ & $n$ & $m^\mathrm{T}$   & $\mbox{GPP}_{\rm m}$ & \eqref{RSfix}  & $\mbox{GPP}_{\rm QAP}$  & $\mbox{GPP}_{\rm m-ind}$  & GPP$_{\rm fix}$   \\
\noalign{\smallskip}\hline\noalign{\smallskip}
$J(6,2)$& 15 & (8,7) & 23 & 23  & 23  & 26 &  24  \\
$J(7,2)$& 21 & (12,9) & 36 & 37  & 36  & 38 &  38  \\
$J(9,2)$& 36 & (26,10) & 65 & 66 & 65  & 65 &  67  \\
{\tt HS} & 50 & (46,4) & 19  & 19 & 19  & 19  & 21 \\
{\tt Gewirtz} & 56 & (53,3) & 23  &23  & 24  &  23 & 26  \\
$J(12,2)$& 66 & (33,33) & 198 & 199& 198  & 198 & 199 \\
$M_{22}$ & 77 & (74,3)& 41 & 41  & 42 & 41  & 44   \\
$J(15,2)$& 105 & (85,20) & 243 & 243 & 243  & 243  & 246 \\
\noalign{\smallskip}\hline
\end{tabular}
\end{center}
\end{table}

\begin{table}[h!]
\caption{Computational times (s.) for the min GBP.}  \label{tab:GBPtime}
\begin{center}
\begin{tabular}{ccccccc}
\hline\noalign{\smallskip}
$G$ & $n$    & \eqref{RSfix}  & $\mbox{GPP}_{\rm QAP}$   & $\mbox{GPP}_{\rm m-ind}$  &  GPP$_{\rm fix}$ \\
\noalign{\smallskip}\hline\noalign{\smallskip}
$J(6,2)$& 15  & 0.42 & 0.20  & 0.17  & 2.30 \\
$J(7,2)$& 21  & 0.66 & 0.23  & 0.18  & 2.70 \\
$J(9,2)$& 36  & 1.06 & 0.48  & 0.23  & 5.21 \\
HS      & 50  & 0.68 & 0.62  & 0.30  & 6.58\\
Gewirtz & 56  & 1.67 & 1.66  & 0.34 & 15.87 \\
$J(12,2)$& 66 & 1.02 & 0.95  & 0.34  & 8.87\\
$M_{22}$ & 77 & 1.55 & 3.15  & 0.30 & 19.12\\
$J(15,2)$& 105& 2.27 & 3.41  & 0.54 & 29.19\\
\noalign{\smallskip}\hline
\end{tabular}
\end{center}
\end{table}

\subsection{Aggregated triangle and independent set constraints} \label{sec:Aggregate}

In Section \ref{sec:aggregateConstr} we showed how to aggregate triangle and independent set constraints for
the case that the data matrices of $\mbox{GPP}_{\rm m}$ belong to a coherent algebra, and that this is
efficient when the rank of this algebra is small. In this section we provide numerical results for  graphs
whose adjacency matrices indeed belong to a coherent algebra of small rank. In particular, besides
the Johnson graph $J(7,2)$, we consider the distance-regular {\tt Pappus}, {\tt Desargues}, {\tt Foster}, and
{\tt Biggs-Smith} graphs (see \cite{bcn}), as well as the {\tt Dyck} graph (the graph on the triangles of the
Shrikhande graph, where two triangles are adjacent if they share an edge). In the third column of Table
\ref{tab:Aggregate}, we list the rank of the smallest coherent configuration containing the corresponding
adjacency matrix (i.e., coming from the combinatorial symmetry). In all cases, this coherent configuration is
the same as the one coming from the automorphism group. In columns five to eight, we list bounds obtained by
solving $\mbox{GPP}_{\rm m}$, $\mbox{GPP}_{\rm m\triangle}$, $\mbox{GPP}_{\rm m-ind}$, and $\mbox{GPP}_{\rm
m}$ with (aggregated) triangle and independent set inequalities, respectively.

The numerical results in
\cite{Sot11} show that to solve $\mbox{GPP}_{\rm m\triangle-ind}$ for a graph without symmetry and
 $n=100$ takes more than 3 hours. However, each bound presented in Table \ref{tab:Aggregate}  is computed in {\em less} than a second (!).

The results here also show that, for most of the cases, adding triangle inequalities to $\mbox{GPP}_{\rm m}$
increases the bound {\em more} than adding the independent set inequalities to $\mbox{GPP}_{\rm m}$. Note that
$J(7,2)$ is a strongly regular graph for which $\mbox{GPP}_{\rm m}$ improves after adding all independent set
constraints. \\

\begin{table}[h!]
\caption{Lower bounds for the min GBP. }  \label{tab:Aggregate}
\begin{center}
\begin{tabular}{cccccccc}
\noalign{\smallskip}\hline\noalign{\smallskip}
$G$ & $n$ & $r_{\rm aut}$  & $m^\mathrm{T}$ &  $\mbox{GPP}_{\rm m}$  & $\mbox{GPP}_{\rm m\triangle}$   & $\mbox{GPP}_{\rm m-ind}$ & $\mbox{GPP}_{\rm m\triangle-ind}$ \\
\noalign{\smallskip}\hline\noalign{\smallskip}
{\tt Pappus} & 18 & 5 &  $(10,8)$ & 6 & 7 & 7  &  7 \\
{\tt Desargues} & 20 & 6 & $(15,5)$ & 4  & 5  &  4 & 5 \\
$J(7,2)$   &  21 & 3 & $(11,10)$ & 37 & 37 & 40  & 40 \\
{\tt Dyck} & 32 & 10 & $(16,16)$  & 7 & 8 & 7 & 8 \\
{\tt Foster}    & 90  & 9 & $(45, 45)$ & 13 & 18 & 14 &  19 \\
{\tt Biggs-Smith} & 102 & 8 & $(70,32)$ &  10 & 15 & 10 & 15 \\
\noalign{\smallskip}\hline
\end{tabular}
\end{center}
\end{table}

In Table \ref{tab:Aggregate22}, we also list the eigenvalue bound, \eqref{RSfix}, $\mbox{GPP}_{\rm QAP}$, and
GPP$_{\rm fix}$ for the same problems as in Table \ref{tab:Aggregate}. Due to memory restrictions, we
couldn't compute GPP$_{\rm fix}$ for the {\tt Foster} and {\tt Biggs-Smith} graph. For these graphs, we
computed \eqref{RSfix} without exploiting their symmetry.  \\

The computational results show that in all cases the eigenvalue bound coincides with $\mbox{GPP}_{\rm m}$.
The results also show that $\mbox{GPP}_{\rm QAP}$ equals
$\mbox{GPP}_{\rm m}$ in all cases except for the {\tt Desargues} graph, and that for the listed graphs fixing
edges does not improve $\mbox{GPP}_{\rm m}$ and/or $\mbox{GPP}_{\rm QAP}$ while  adding triangle and/or
independent set constraints does.

\begin{table}[h!]
\caption{Lower bounds for the min GBP. }  \label{tab:Aggregate22}
\begin{center}
\begin{tabular}{cccccc}
\noalign{\smallskip}\hline\noalign{\smallskip}
$G$ & $m^\mathrm{T}$ & eig & \eqref{RSfix} & $\mbox{GPP}_{\rm QAP}$  &  GPP$_{\rm fix}$  \\[1ex]
\noalign{\smallskip}\hline\noalign{\smallskip}
{\tt Pappus}    & $(10,8)$  & 6  & 6 & 6  & 6   \\
{\tt Desargues} & $(15,5)$  & 4  & 4 & 5  & 6  \\
$J(7,2)$        & $(11,10)$ & 37 & 38 & 37 & 38  \\
{\tt Dyck}      & $(16,16)$ & 7  & 7 & 7  & 7  \\
{\tt Foster}    & $(45, 45)$& 13 & 13 & 13 &  --  \\
{\tt Biggs-Smith}& $(70,32)$ & 10 & 10 & 10 &  -- \\
\noalign{\smallskip}\hline
\end{tabular}
\end{center}
\end{table}

\section{Conclusion}

In this paper, we presented several new bounds for the graph partition problem and also showed how to simplify existing ones, when possible.

In particular, in Theorem \ref{thm:AnyGraph} we derived an eigenvalue bound for the GPP that is applicable for any graph, extending a well-known result in spectral graph theory.
Further, we simplified the relaxation ${\rm GPP}_{\rm m}$
for different classes of graphs such as strongly regular graphs and certain Johnson and Kneser graphs.
We showed how to strengthen ${\rm GPP}_{\rm m}$ by aggregating triangle and independent set constraints when possible,
which leads to huge improvements in computational abilities and times.
The numerical results show that, in general, adding triangle inequalities to ${\rm GPP}_{\rm m}$ strengthened the bound more
than adding the independent set constraints.

Finally, we showed how to strengthen the matrix and vector lifting relaxation for the GPP,
 i.e., ${\rm GPP}_{\rm m}$ and ${\rm GPP}_{\rm QAP}$, respectively,  by adding constraints
that correspond to assigning two vertices of the graph to different parts of the partition.
Such an approach performs very well on highly symmetric graphs when other SDP and eigenvalue-based relaxations provide trivial or weak bounds.

\vspace{1cm}
\noindent
{\bf Acknowledgments}\\
The authors would like to thank two anonymous referees for suggestions that led to an improvement of this paper.

\end{document}